\documentclass{amsart}
\usepackage{amsfonts,amsmath,amsthm,amssymb}
\newtheorem{theorem}{Theorem}

\newtheorem{corollary}{Corollary}
\newtheorem{proposition}{Proposition}
\title[Integral group ring of the Mathieu simple group]
      {Integral group ring \\ of the Mathieu simple group $M_{23}$}

\author[Bovdi]{V.A.~Bovdi}

\author[Konovalov]{A.B.~Konovalov}

\dedicatory{Institute of Mathematics, University of Debrecen, Hungary \\
School of Computer Science, University of St Andrews, Scotland}

\thanks{The research was supported by OTKA grants No.T 037202,
No.T 038059 and Francqui Stichting (Belgium) grant ADSI107. \\
V.A.~Bovdi: Institute of Mathematics, University of Debrecen,
P.O.  Box 12, H-4010 Debrecen, Hungary;
Institute of Mathematics and Informatics, College of Ny\'\i regyh\'aza,
S\'ost\'oi \'ut 31/b, H-4410 Ny\'\i regyh\'aza, Hungary;
E-mail: vbovdi@math.klte.hu. \\
A.B.~Konovalov:
School of Computer Science, University of St Andrews,
Jack Cole Building, North Haugh, St Andrews, Fife, KY16 9SX, Scotland;
E-mail: konovalov@member.ams.org}



\begin{document}

\begin{abstract}
We investigate the classical Zassenhaus conjecture for the unit
group of the integral group ring of Mathieu simple group $M_{23}$
using the Luthar-Passi method. This work is a continuation of the
research that we carried out for Mathieu groups $M_{11}$ and
$M_{12}$. As a consequence, for this group we confirm Kimmerle's
conjecture on prime graphs.
\end{abstract}

\maketitle


\section{Introduction  and main results}

Let $V(\mathbb Z G)$ be the normalized unit group of the integral
group ring $\mathbb Z G$ of a finite group $G$. With respect to
the structure of $V(\mathbb Z G)$, H.~Zassenhaus in
\cite{Zassenhaus} proposed the following conjecture:

\begin{itemize}
\item[] {\bf (ZC)} \qquad Every torsion unit $u \in V(\mathbb ZG)$
is conjugate within the rational group algebra $\mathbb QG$ to an
element of $G$.
\end{itemize}

In the general case the problem is still open, despite a lot of
known results about {\bf (ZC)} and methods developed for its
investigation. In the present paper we will be focused on those of
them which are relevant to finite simple groups.

The first result on {\bf (ZC)} for simple groups was obtained by
I.S.~Luthar and I.B.S.~Passi in \cite{Luthar-Passi} for the
alternating group $A_5$, and the technique used in this paper is
now known as the Luthar-Passi method. Later M.~Hertweck in
\cite{Hertweck1} extended this method and applied it for the
investigation of {\bf (ZC)} for $PSL(2,p^{n})$. This method proved
to be useful for non-simple groups as well. For some recent
results we refer to
\cite{Bovdi-Jespers-Konovalov,Bovdi-Hofert-Kimmerle,
Bovdi-Konovalov,Hertweck1, Hertweck2, Hertweck3, Hofert-Kimmerle}.
Also some related properties can be found in
\cite{Artamonov-Bovdi,Luthar-Trama} and
\cite{Bleher-Kimmerle,Kimmerle}. In the latter papers some
weakened variations of {\bf (ZC)} have been made.

In \cite{Kimmerle} W.~Kimmerle proposed to relate {\bf (ZC)} with
some properties of graphs associated with groups. Let $\# (G)$ be
the set of all primes dividing the order of $G$. The
Gruenberg-Kegel graph (or the prime graph) of $G$ is the graph
$\pi (G)$ with vertices labeled by the primes in $\# (G)$ and with
an edge from $p$ to $q$ if there is an element of order $pq$ in
the group $G$. The following conjecture was introduced by W.~Kimmerle:

\begin{itemize}
\item[] {\bf (KC)} \qquad If $G$ is a finite group then $\pi(G)
=\pi (V(\mathbb Z G))$.
\end{itemize}

\noindent In \cite{Kimmerle} {\bf (KC)} was solved for finite
solvable groups and finite Frobenius groups. Note that with
respect to the Zassenhaus conjecture the investigation of
Frobenius groups was completed by M.~Hertweck and the first author
in \cite{Bovdi-Hertweck}. In our previous papers
\cite{Bovdi-Jespers-Konovalov, Bovdi-Konovalov,
Bovdi-Konovalov-Siciliano} we confirmed {\bf (KC)} for sporadic
simple groups $M_{11}$, $M_{12}$ and some Janko simple groups.

In the present paper we continue these investigations for the
Mathieu simple group $M_{23}$ using the Luthar-Passi method. In
our main result we obtain a lot of information on possible torsion
units in $V(\mathbb Z M_{23})$, and as a consequence of part (i)
of the main theorem we solve {\bf (KC)} for $M_{23}$.

First we introduce some notation. Put $G=M_{23}$. Let
\[
\begin{split}
\mathcal{C} =\{
C_{1}, C_{2a},C_{3a}, & C_{4a}, C_{5a},C_{6a},C_{7a},C_{7b},C_{8a},\\
& C_{11a},C_{11b},C_{14a},C_{14b},
C_{15a},C_{15b},C_{23a},C_{23b}\},
\end{split}
\]
be the collection of all conjugacy classes of $G$, where the first
index denotes the order of the elements of this conjugacy class
and $C_{1}=\{ 1\}$. Suppose $u=\sum \alpha_g g \in V(\mathbb Z G)$
has finite order $k$. Denote by
$\nu_{nt}=\nu_{nt}(u)=\varepsilon_{C_{nt}}(u)=\sum_{g\in C_{nt}}
\alpha_{g}$ the partial augmentation of $u$ with respect to
$C_{nt}$. From S.D.~Berman's Theorem \cite{Berman} one knows that
$\nu_1=\alpha_{1}=0$ and
\begin{equation}
\label{E:1} \sum_{C_{nt}\in \mathcal{C}} \nu_{nt}=1.
\end{equation}
Hence, for any character $\chi$ of $G$, we get that
$\chi(u)=\sum\nu_{nt}\chi(h_{nt})$, where $h_{nt}$ is a
representative of a conjugacy class $ C_{nt}$.

The main result is the following.

\begin{theorem}\label{T:1}
Let $V(\mathbb ZG)$ be the normalized unit group of the integral
group ring $\mathbb ZG$, where $G$ is the sporadic simple Mathieu
group $M_{23}$. Let $u$ be  a torsion unit of $V(\mathbb ZG)$ of
order $|u|$. Then:

\begin{itemize}

\item[(i)] if $|u| \not \in \{12,24,28,56\}$, then $|u|$ coincides
with the order of some element $g\in G$. Equivalently, there is no
elements of orders $10$, $21$, $22$, $33$, $35$, $46$, $55$, $69$,
$77$, $115$, $161$ and $253$ in $V(\mathbb ZG)$.

\item[(ii)] if $|u| \in \{2,3,5,23\}$, then $u$ is rationally
conjugate to some  $g\in G$;

\item[(iii)] if $|u|=4$, then the tuple of partial augmentations
of $u$ belongs to the set
\[
\begin{split}
\big \{\;  (\nu_{2a}, \nu_{3a},\;& \nu_{4a}, \nu_{5a}, \nu_{6a},
\nu_{7a},\nu_{7b}, \nu_{8a}, \nu_{11a}, \nu_{11b},  \nu_{14a},
\nu_{14b}, \nu_{15a},\\
& \nu_{15b}, \nu_{23a}, \nu_{23b} ) \in \mathbb Z^{16}
\quad \mid \quad \nu_{kx}=0, \\
&kx \not \in \{2a, 4a\},  \quad (\nu_{2a}, \nu_{4a} ) \in \{\; (
0, 1 ), \; ( -2, 3 ), \; ( 2, -1 ) \; \} \; \big \};
\end{split}
\]

\item[(iv)] if $|u|=6$,  then  the tuple of partial augmentations
of $u$ belongs to the set
\[
\begin{split}
\big \{\; & (\nu_{2a}, \nu_{3a}, \nu_{4a}, \nu_{5a}, \nu_{6a},
\nu_{7a},\nu_{7b}, \nu_{8a}, \nu_{11a}, \nu_{11b},  \nu_{14a},
\nu_{14b}, \nu_{15a}, \\
& \nu_{15b}, \nu_{23a}, \nu_{23b} ) \in \mathbb Z^{16}
\quad \mid \nu_{kx}=0,\quad   kx\not\in \{2a, 3a, 6a\}, \\
& (\nu_{2a}, \nu_{3a}, \nu_{6a}) \in \{\;    ( 0, 0, 1 ),\quad  ( 4, -6, 3 ),\quad  ( -6, 6, 1 ),\quad  ( -6, 9, -2 ), \\
& ( -6, 12, -5 ),\quad  ( -4, 3, 2 ),\quad  ( -4, 6, -1 ),\quad 
( -4, 9, -4 ),\quad  ( -2, 0, 3 ),\quad \\
& ( -2, 3, 0 ), \qquad ( -2, 6, -3 ),\quad  ( 0, -3, 4 ),\quad  (
0, 3, -2 ),\qquad  ( 2, -6, 5 ), \\
& ( 2, -3, 2 ),\qquad  ( 2, 0, -1 ),\qquad (4,-9, 6 ), \quad  ( 4,
-3, 0 ),\qquad  ( 6, -12, 7 ), \\
& ( 6, -9, 4 ),\qquad    ( 6, -6, 1 )  \; \} \qquad  \big  \};
\end{split}
\]

\item[(v)] if $|u|=7$, then the tuple of partial augmentations of
$u$ belongs to the set
\[
\begin{split}
\big \{\quad & (\nu_{2a}, \nu_{3a}, \nu_{4a}, \nu_{5a}, \nu_{6a},
\nu_{7a},\nu_{7b}, \nu_{8a}, \nu_{11a}, \nu_{11b}, \nu_{14a},
\nu_{14b}, \nu_{15a}, \\
& \nu_{15b}, \nu_{23a}, \nu_{23b} ) \in \mathbb Z^{16}
\quad \mid \nu_{kx}=0,\\
& kx \not \in \{7a, 7b\},\quad (\nu_{7a}, \nu_{7b} ) \in \{\; (0,1
), (2,-1), (1,0),\;  (-1,2) \; \}\quad \big \};
\end{split}
\]

\item[(vi)] if $|u|=8$, then the tuple of partial augmentations of
$u$ belongs to the set
\[
\begin{split}
\big \{\quad (\nu_{2a}, \nu_{3a}, & \nu_{4a}, \nu_{5a}, \nu_{6a},
\nu_{7a}, \nu_{7b}, \nu_{8a}, \nu_{11a}, \nu_{11b}, \nu_{14a},
\nu_{14b}, \nu_{15a}, \\
&  \nu_{15b}, \nu_{23a}, \nu_{23b} ) \in \mathbb Z^{16}
\quad \mid \nu_{kx}=0,\\
&kx\not\in \{2a, 4a, 8a\},\quad (\nu_{2a}, \nu_{4a}, \nu_{8a} )
\in
\{\quad  (-2,4,-1), \\
& ( 0, 4, -3 ), \quad ( 2, 0, -1 ), \quad
  ( 0, 0, 1 ),\quad  ( 0, 2, -1 ),\quad
  (2,-2,1), \\
& ( -2, 2, 1 ),\quad  ( 0, -2, 3 ), \quad (-2,6,-3), \quad
(2,-4,3) \}\quad \big \};
\end{split}
\]

\item[(vii)] if $|u|=11$, then the tuple of partial augmentations
of $u$ belongs to the set
\[
\begin{split}
\big \{ \; (\nu_{2a}, & \nu_{3a}, \nu_{4a}, \nu_{5a}, \nu_{6a},
\nu_{7a},\nu_{7b}, \nu_{8a}, \nu_{11a}, \nu_{11b}, \nu_{14a},
\nu_{14b}, \nu_{15a},\\
& \nu_{15b}, \nu_{23a}, \nu_{23b} ) \in \mathbb Z^{16} \quad \mid\quad \nu_{kx}=0,\\
& kx\not\in \{11a, 11b\},\quad  (\nu_{11a}, \nu_{11b})
\in \{\; ( -5, 6 ),\; ( 5, -4 ),\; ( 0, 1 ), \;( -2, 3 ), \\
& ( -8, 9 ),\; ( -9, 10 ),\; ( -6, 7 ),\; ( 2, -1 ),\; ( 8, -7
),\; ( -3, 4 ),\;( 6, -5 ),\; ( -4, 5 ), \\
& ( 1, 0 ),\; ( 3, -2 ),\; ( 9, -8 ),\;  ( 10, -9 ),\;  ( 7, -6 ),
\; ( -1, 2 ),\;  ( -7, 8 ), \;( 4, -3 ) \; \}\quad \big  \};
\end{split}
\]

\item[(viii)] if $|u|=15$, then the tuple of partial augmentations
of $u$ belongs to the set
\[
\begin{split}
\big \{ \; (\nu_{2a}, & \nu_{3a}, \nu_{4a}, \nu_{5a}, \nu_{6a},
\nu_{7a}, \nu_{7b}, \nu_{8a}, \nu_{11a}, \nu_{11b},  \nu_{14a},
\nu_{14b}, \nu_{15a},\\
&
\nu_{15b}, \nu_{23a}, \nu_{23b} ) \in \mathbb Z^{16} \quad \mid\quad \nu_{kx}=0,\\
& kx\not\in \{3a, 5a, 15a, 15b\},\quad (\nu_{3a}, \nu_{5a},
\nu_{15a}, \nu_{15b} )\in \{ \;
( -3, 5, -1, 0 ),\\
&( -3, 5, 0, -1 ),\;( 0, 0, 0, 1 ),\; ( 0, 0, 1, 0 ),\;
  ( 3, -5, 1, 2 ), \;( 3, -5, 2, 1 )
 \;\}\quad \big  \}.
\end{split}
\]

\end{itemize}

\end{theorem}

\noindent As an immediate consequence of the first part of the
Theorem one obtains Kimmerle's conjecture for $M_{23}$.

\begin{corollary} If $G=M_{23}$ then
$\pi(G)=\pi(V(\mathbb ZG))$.
\end{corollary}


\section{Preliminaries}

Throughout the paper we simply denote $M_{23}$ by $G$. The
$p$-Brauer character table of the group $G$ will be denoted by
$\mathfrak{BCT}{(p)}$.

The crucial restriction on partial augmentations is given by the
next result, detailed explanation of which can be found in
\cite{Luthar-Passi} and \cite{Hertweck1}.

\begin{proposition}\label{P:1}
(see \cite{Luthar-Passi,Hertweck1}) Let either $p=0$ or $p$ is a
prime divisor of $|G|$. Suppose that $u\in V( \mathbb Z G) $ has
finite order $k$ and assume $k$ and $p$ are coprime in case $p\neq
0$. If $z$ is a primitive $k$-th root of unity and $\chi$ is
either a classical character or a $p$-Brauer character of $G$
then, for every integer $l$, the number
\begin{equation}
\label{E:2} \mu_l(u,\chi, p )=\frac{1}{k} \sum_{d|k}Tr_{  \mathbb
Q (z^d)/  \mathbb Q } \{\chi(u^d)z^{-dl}\}
\end{equation}
is a non-negative integer.
\end{proposition}

\noindent If $p=0$, we will use the notation $\mu_l(u,\chi, * )$
for $\mu_l(u,\chi , 0)$.

For the determination of possible orders of units we use the next
bound.

\begin{proposition}\label{P:2}  (see  \cite{Cohn-Livingstone})
The order of a torsion element $u\in V(\mathbb ZG)$ is a divisor
of the exponent of $G$.
\end{proposition}

\noindent The next two results, the first of which is a special
case of the second one, already yield that several partial
augmentations are zero.

\begin{proposition}\label{P:3}(see \cite{Luthar-Passi} and
Theorem 2.7 in \cite{Marciniak-Ritter-Sehgal-Weiss}) Let $u$ be a
torsion unit of $V(\mathbb ZG)$. Let $C$ be a conjugacy class of
$G$. If $a\in C$ and $p$ is a prime dividing  the order of $a$ but
not the order of $u$ then  $\varepsilon_C(u)=0$.
\end{proposition}

\begin{proposition}\label{P:4}
(see \cite{Hertweck2}, Proposition 3.1; \cite{Hertweck1}, Proposition 2.2)
Let $G$ be a finite group and let $u$ be a torsion unit in
$V(\mathbb ZG)$. If $x$ is an element of $G$ whose $p$-part, for
some prime $p$, has order strictly greater than the order of the
$p$-part of $u$, then $\varepsilon_x(u)=0$.
\end{proposition}

\noindent The following result relates the solution of the
Zassenhaus conjecture to partial augmentations of torsion units.

\begin{proposition}\label{P:5}
(see \cite{Luthar-Passi} and Theorem 2.5 in
\cite{Marciniak-Ritter-Sehgal-Weiss}) Let $u\in V(\mathbb Z G)$ be
of order $k$. Then $u$ is conjugate in $\mathbb QG$ to an element
$g \in G$ if and only if for each $d$ dividing $k$ there is
precisely one conjugacy class $C$ with partial augmentation
$\varepsilon_{C}(u^d) \neq 0 $.
\end{proposition}

\noindent Finally, the next result is useful for the investigation
of $p$-elements in $V(\mathbb ZG)$.

\begin{proposition}\label{P:6}
(see \cite{Cohn-Livingstone}) Let $p$ be a prime, and let $u$ be a
torsion unit of $V(\mathbb ZG)$ of order $p^n$. Then for $m \ne n$
the sum of all partial augmentations of $u$ with respect to
conjugacy classes of elements of order $p^m$ is divisible by $p$.
\end{proposition}


\section{Proof of the Theorem}

It is well known \cite{GAP,Gorenstein} that $|G|=2^7 \cdot 3^2
\cdot 5 \cdot 7 \cdot 11 \cdot 23$ and  $exp(G) = 2^3 \cdot 3
\cdot 5 \cdot 7 \cdot 11 \cdot 23$.
The  character table of $G$, as well as the Brauer character
tables $\mathfrak{BCT}{(p)}$, where $p \in \{2,3,5,7,11,23\}$, can
be found using the GAP system \cite{GAP}. Throughout the paper we
will use the notation, indexation inclusive, for the characters
and conjugacy classes as used in GAP.

Since the group $G$ possesses elements of orders $2$, $3$, $4$,
$5$, $6$, $7$, $8$, $11$, $14$, $15$ and $23$, first of
all we shall investigate units of these orders. After this, by
Proposition \ref{P:2}, the order of each torsion unit divides the
exponent of $G$, so it will be enough to consider units of orders
$10$, $12$, $21$, $22$, $24$, $28$, $33$, $35$, $46$, $55$, $56$, 
$69$, $77$, $115$, $161$ and  $253$, because if $u$ will be a unit 
of another possible order, then there is $t \in \mathbb N$ such that $u^t$
has an order from this list. For all these orders except $12$, $24$, $28$
and $56$ we will prove that units of such order do not appear in 
$V(\mathbb ZG)$.

To reduce the volume of the paper, we omit from consideration
units of orders 12, 14 and 28, because in these cases we obtain
systems with enormous number of solutions. Note that with the
LAGUNA package \cite{LAGUNA} we can investigate units of these orders and show
that there are no units of order 28, and, consequently, of order 56 in $V(\mathbb ZG)$.
Unfortunately, for units of order 12 the Luthar-Passi method is
not enough to prove the same result.

Now we consider each case separately:

\noindent $\bullet$ Let $u$ be a unit of order $2$, $3$ or $5$.
Then using Propositions \ref{P:3} and \ref{P:4} we obtain that all
partial augmentations except one are zero.

\noindent $\bullet$ Let $u$ be a unit of order $4$. Then \;
$\nu_{2a}+\nu_{4a}=1$ \; by \eqref{E:1} and Proposition~\ref{P:4}.
By \eqref{E:2}, using $\mathfrak{BCT}{(23)}$ we obtain the system
of inequalities
\[
\begin{split}
\mu_0(u,\chi_{10},23) & = \textstyle \frac{1}{4} (-16 \nu_{2a} +
272) \geq 0; \quad
\mu_2(u,\chi_{10},23)   = \textstyle \frac{1}{4} (16 \nu_{2a} + 272) \geq 0; \\
\mu_0(u,\chi_2,23) & =\textstyle \frac{1}{4} ( 10 \nu_{2a} + 2 \nu_{4a} + 26 ) \geq 0;\\
\mu_2(u,\chi_2,23) & =\textstyle \frac{1}{4} ( -10 \nu_{2a}  -2 \nu_{4a} + 26 ) \geq 0;\\
\mu_0(u,\chi_3,23) & = \textstyle \frac{1}{4} ( -6 \nu_{2a} + 2 \nu_{4a} + 42 ) \geq 0;\\
\mu_2(u,\chi_3,23) & = \textstyle \frac{1}{4} ( 6 \nu_{2a} -2 \nu_{4a} + 42 ) \geq 0,\\
\end{split}
\]
that has only three integral solutions $(\nu_{2a},\nu_{4a}) \in
\{\; (0,1),\; (-2,3),\; (2,-1)\; \}$ satisfying Proposition
\ref{P:6}, such that all $\mu_i(u,\chi_j,23)$ are non-negative
integers.

\noindent$\bullet$ Let $u$ be a unit of order $6$. Clearly,
$\nu_{2a}+\nu_{3a}+\nu_{6a}=1$ by \eqref{E:1} and Proposition~\ref{P:4}. 
Now from the system of inequalities
\[
\begin{split}
\mu_{0}(u,\chi_{3},*) & = \textstyle \frac{1}{6} (-6\nu_{2a} + 42)
\geq 0; \qquad
\mu_{3}(u,\chi_{3},*) = \textstyle \frac{1}{6} (6\nu_{2a} + 48) \geq 0; \\
\mu_{0}(u,\chi_{12},*) & = \textstyle \frac{1}{6} (- 8 \nu_{3a} +
888) \geq 0; \quad \;
\mu_{3}(u,\chi_{12},*) = \textstyle \frac{1}{6} (8 \nu_{3a} + 888) \geq 0; \\
\mu_{1}(u,\chi_{2},*) & = \textstyle \frac{1}{6} (6\nu_{2a} + 4 \nu_{3a} + 12) \geq 0; \\
\mu_{2}(u,\chi_{2},*) & = \textstyle \frac{1}{6} (-6\nu_{2a} - 4 \nu_{3a} + 24) \geq 0; \\
\mu_{1}(u,\chi_{5},*) & = \textstyle \frac{1}{6} (22\nu_{2a} + 5 \nu_{3a} +  \nu_{6a} + 203) \geq 0; \\
\mu_{0}(u,\chi_{4},7) & = \textstyle \frac{1}{6} (32\nu_{2a} + 2 \nu_{3a} + 2 \nu_{6a} + 226) \geq 0; \\
\mu_{3}(u,\chi_{2},23) & = \textstyle \frac{1}{6} (-10\nu_{2a} - 6 \nu_{3a} + 2 \nu_{6a} + 22) \geq 0, \\
\end{split}
\]
we obtain only those integral solutions that are listed in part
(iv) of the Theorem.

\noindent $\bullet$ Let $u$ be a unit of order $7$. By \eqref{E:1}
and Proposition \ref{P:4} we have $\nu_{7a}+\nu_{7b}=1$. Now using
\eqref{E:2} we obtain the system of inequalities
\[
\begin{split}
\mu_{1}(u,\chi_{3},*) & = \textstyle \frac{1}{7} (4\nu_{7a} - 3 \nu_{7b} + 45) \geq 0; \\
\mu_{3}(u,\chi_{3},*) & = \textstyle \frac{1}{7} (-3\nu_{7a} + 4 \nu_{7b} + 45) \geq 0; \\
\mu_{1}(u,\chi_{2},2) & = \textstyle \frac{1}{7} (-4\nu_{7a} + 3 \nu_{7b} + 11) \geq 0; \\
\mu_{3}(u,\chi_{2},2) & = \textstyle \frac{1}{7} (3\nu_{7a} - 4 \nu_{7b} + 11) \geq 0,\\
\end{split}
\]
that has solutions $(\nu_{7a}, \nu_{7b}) \in \{\; ( 0, 1 ),\; ( 2,
-1 ),\; ( 1,0),\; (-1, 2 )  \; \}. $

\noindent $\bullet$ Let $u$ be a unit of order $8$. We get
$\nu_{2a}+\nu_{4a}+\nu_{8a}=1$ by \eqref{E:1} and Proposition~\ref{P:4}. 
It is necessary to consider three cases, defined by
part (iii) of the Theorem. By \eqref{E:2} we obtain the system of
inequalities
\[
\begin{split}
\mu_{0}(u,\chi_{4},7) & = \textstyle \frac{1}{8} (64\nu_{2a} +
\alpha_1) \geq 0; \quad
\mu_{4}(u,\chi_{4},7)   = \textstyle \frac{1}{8} (-64\nu_{2a} + \alpha_1) \geq 0; \\
\mu_{0}(u,\chi_{2},*) & = \textstyle \frac{1}{8} (24\nu_{2a} + 8 \nu_{4a} + \alpha_2) \geq 0; \\
\mu_{4}(u,\chi_{2},*) & = \textstyle \frac{1}{8} (-24\nu_{2a} - 8 \nu_{4a} + \alpha_2) \geq 0; \\
\mu_{0}(u,\chi_{3},*) & = \textstyle \frac{1}{8} (-12\nu_{2a} + 4 \nu_{4a} - 4 \nu_{8a} + \alpha_3) \geq 0; \\
\mu_{4}(u,\chi_{3},*) & = \textstyle \frac{1}{8} (12\nu_{2a} - 4 \nu_{4a} + 4 \nu_{8a} + \alpha_3) \geq 0; \\
\mu_{0}(u,\chi_{2},23) & = \textstyle \frac{1}{8} (20\nu_{2a} + 4 \nu_{4a} - 4 \nu_{8a} + \alpha_4) \geq 0; \\
\mu_{4}(u,\chi_{2},23) & = \textstyle \frac{1}{8} (-20\nu_{2a} - 4 \nu_{4a} + 4 \nu_{8a} + \alpha_4) \geq 0, \\
\end{split}
\]
where $(\alpha_1,\alpha_2,\alpha_3,\alpha_4) =
\begin{cases}
(224,32,44,28), & \text{if} \quad \chi(u^2)=\chi(4a) ;\\
(160,16,60,12), & \text{if} \quad \chi(u^2)=-2\chi(2a)+3\chi(4a) ;\\
(288,48,28,44), & \text{if} \quad \chi(u^2)=2\chi(2a)-\chi(4a) .\\
\end{cases}
$ \vspace{10pt}

\noindent Union of solutions for all three cases, taking into
account restrictions from Proposition \ref{P:6}, gives us part
(vi) of the Theorem.

\noindent $\bullet$ Let $u$ be a unit of order $11$. Clearly,
$\nu_{11a}+\nu_{11b}=1$ by \eqref{E:1} and Proposition~\ref{P:4}.
Using $\mathfrak{BCT}{(3)}$ by \eqref{E:2} we obtain the system of
inequalities
\[
\begin{split}
\mu_{1}(u,\chi_{5},3) & = \textstyle \frac{1}{11} (6 \nu_{11a} - 5 \nu_{11b} + 104) \geq 0; \\
\mu_{2}(u,\chi_{5},3) & = \textstyle \frac{1}{11} (-5 \nu_{11a} + 6 \nu_{11b} + 104) \geq 0, \\
\end{split}
\]
that has only 20 integer solutions listed in part (vi) of the
Theorem, such that are $\mu_{1}(u,\chi_{5},3)$ and
$\mu_{2}(u,\chi_{5},3)$ are non-negative integers.

\noindent $\bullet$ Let $u$ be a unit of order $15$. Again,
$\nu_{3a}+\nu_{5a}+\nu_{15a}+\nu_{15b}=1$ by \eqref{E:1} and
Proposition~\ref{P:4}. Using \eqref{E:2}, we obtain the system of
inequalities
\[
\begin{split}
\mu_{0}(u,\chi_{7},2)   & = \textstyle \frac{1}{15} (-40 \nu_{3a}
+ 210) \geq 0; \quad
\mu_{5}(u,\chi_{7},2)     = \textstyle \frac{1}{15} (20 \nu_{3a} + 225) \geq 0; \\
\mu_{0}(u,\chi_{9},2)   & = \textstyle \frac{1}{15} (- 24 \nu_{5a}
+ 240) \geq 0; \quad
\mu_{5}(u,\chi_{9},2)     = \textstyle \frac{1}{15} (+ 12 \nu_{5a} + 240) \geq 0; \\
\mu_{0}(u,\chi_{5},*) & = \textstyle \frac{1}{15} (40\nu_{3a} + 240) \geq 0; \\
\end{split}
\]

\[
\begin{split}
\mu_{0}(u,\chi_{2},*) & = \textstyle \frac{1}{15} (32\nu_{3a} + 16 \nu_{5a} - 8 \nu_{15a} - 8 \nu_{15b} + 38) \geq 0; \\
\mu_{3}(u,\chi_{2},*) & = \textstyle \frac{1}{15} (-8\nu_{3a} - 4 \nu_{5a} + 2 \nu_{15a} + 2 \nu_{15b} + 28) \geq 0; \\
\mu_{5}(u,\chi_{2},*) & = \textstyle \frac{1}{15} (-16\nu_{3a} - 8 \nu_{5a} + 4 \nu_{15a} + 4 \nu_{15b} + 26) \geq 0; \\
\mu_{1}(u,\chi_{7},*) & = \textstyle \frac{1}{15} (-3\nu_{3a} +  \nu_{5a} + 7 \nu_{15a} - 8 \nu_{15b} + 233) \geq 0; \\
\mu_{7}(u,\chi_{7},*) & = \textstyle \frac{1}{15} (-3\nu_{3a} +  \nu_{5a} - 8 \nu_{15a} + 7 \nu_{15b} + 233) \geq 0; \\
\mu_{1}(u,\chi_{2},2) & = \textstyle \frac{1}{15} (2\nu_{3a} +  \nu_{5a} + 7 \nu_{15a} - 8 \nu_{15b} + 8) \geq 0; \\
\mu_{7}(u,\chi_{2},2) & = \textstyle \frac{1}{15} (2\nu_{3a} +  \nu_{5a} - 8 \nu_{15a} + 7 \nu_{15b} + 8) \geq 0; \\
\mu_{0}(u,\chi_{4},7) & = \textstyle \frac{1}{15} (8\nu_{3a} - 16 \nu_{5a} + 8 \nu_{15a} + 8 \nu_{15b} + 202) \geq 0; \\
\mu_{1}(u,\chi_{10},23) & = \textstyle \frac{1}{15} (\nu_{3a} + 6 \nu_{15a} - 9 \nu_{15b} + 279) \geq 0; \\
\mu_{1}(u,\chi_{12},23) & = \textstyle \frac{1}{15} (-\nu_{3a} - 6 \nu_{15a} + 9 \nu_{15b} + 666) \geq 0. \\
\end{split}
\]
This system has only six solutions such that all
$\mu_{i}(u,\chi_{j},p)$ are non-negative integers, listed in part
(viii) of the Theorem.

\noindent $\bullet$ Let $u$ be a unit of order $23$. By
\eqref{E:1} and Proposition~\ref{P:4}  we have
$\nu_{23a}+\nu_{23b}=1$. Using \eqref{E:2} we obtain the system of
inequalities
\[
\begin{split}
\mu_{1}(u,\chi_{2},2) & = \textstyle \frac{1}{23} (12 \nu_{23a} - 11 \nu_{23b} + 11) \geq 0; \\
\mu_{5}(u,\chi_{2},2) & = \textstyle \frac{1}{23} (-11 \nu_{23a} + 12 \nu_{23b} + 11) \geq 0, \\
\end{split}
\]
that has only two trivial solutions: $(\nu_{23a}, \nu_{23b} )\in
\{ \;( 1, 0 ), ( 0, 1 ) \; \}$.

It remains to prove part (i) of the Theorem.

\noindent $\bullet$ Let $u$ be a unit of order $10$. By
\eqref{E:1} and Proposition~\ref{P:4}  we have
$\nu_{2a}+\nu_{5a}=1$. Using \eqref{E:2} we obtain the system of
inequalities
\[
\begin{split}
\mu_{0}(u,\chi_{5},*) & = \textstyle \frac{1}{10} (88 \nu_{2a} +
252) \geq 0; \quad
\mu_{5}(u,\chi_{5},*)   = \textstyle \frac{1}{10} (-88 \nu_{2a} + 208) \geq 0; \\
\mu_{1}(u,\chi_{5},3) & = \textstyle \frac{1}{10} (8 \nu_{2a} - \nu_{5a} + 97) \geq 0, \\
\end{split}
\]
that has no integral solutions such that all
$\mu_{i}(u,\chi_{j},p)$ are non-negative integers.

\noindent $\bullet$ Let $u$ be a unit of order $21$. We have
$\nu_{3a}+\nu_{7a}+\nu_{7b}=1$ by \eqref{E:1} and Proposition~\ref{P:4}. 
We need to consider four cases given in part (v) of the Theorem:
\[
\begin{matrix}
\text{Case 1.} & \chi(u^3) & = & \chi(7a). \qquad
& \text{Case 3.} & \chi(u^3) & = & 2 \chi(7a) -   \chi(7b). \\
\text{Case 2.} & \chi(u^3) & = & \chi(7b). \qquad &\text{Case 4.}
& \chi(u^3) & = & - \chi(7a) + 2 \chi(7b).
\end{matrix}
\]
In all these cases using \eqref{E:2} we obtain the system of
inequalities
\[
\begin{split}
\mu_{0}(u,\chi_{2},23) & = \textstyle \frac{1}{21} (36 \nu_{3a} +
27) \geq 0; \quad
\mu_{7}(u,\chi_{2},23)   = \textstyle \frac{1}{21} (-18 \nu_{3a} + 18) \geq 0; \\
\mu_{1}(u,\chi_{10},23) & = \textstyle \frac{1}{21} (\nu_{3a} + 279) \geq 0, \\
\end{split}
\]
that has no integral solutions such that all
$\mu_{i}(u,\chi_{j},23)$ are non-negative integers.

\noindent $\bullet$ Let $u$ be a unit of order $22$. We have
$\nu_{2a}+\nu_{11a}+\nu_{11b}=1$ by \eqref{E:1} and Proposition~\ref{P:4}. 
We need to consider 20 cases given in part (vi) of the
Theorem, but in all cases using \eqref{E:2} we obtain the system
of inequalities
$$
\mu_{0}(u,\chi_{2},*) = \textstyle \frac{1}{22} (60 \nu_{2a} + 28)
\geq 0; \quad \mu_{11}(u,\chi_{2},*) = \textstyle \frac{1}{22}
(-60 \nu_{2a} + 16) \geq 0,
$$
that has no integral solutions such that all
$\mu_{i}(u,\chi_{j},*)$ are non-negative integers.

\noindent $\bullet$ Let $u$ be a unit of order $33$. We have
$\nu_{3a}+\nu_{11a}+\nu_{11b}=1$ by \eqref{E:1} and Proposition~\ref{P:4}. 
Again, we have 20 cases, and all of them give us the same system
$$
\mu_{0}(u,\chi_{2},*) = \textstyle \frac{1}{33} (80 \nu_{3a} + 30)
\geq 0; \quad \mu_{11}(u,\chi_{2},*) = \textstyle \frac{1}{33}
(-40 \nu_{3a} + 18) \geq 0,
$$
that has no integral solutions such that all
$\mu_{i}(u,\chi_{j},*)$ are non-negative integers.

\noindent $\bullet$ Let $u$ be a unit of order $35$. We have
$\nu_{5a}+\nu_{11a}+\nu_{11b}=1$ by \eqref{E:1} and Proposition~\ref{P:4}. 
We need to consider 4 cases, given by part (v) of the
Theorem, and all of them give us the same system of inequalities
$$
\mu_{0}(u,\chi_{2},23) = \textstyle \frac{1}{35} (24 \nu_{5a} +
25) \geq 0; \quad \mu_{7}(u,\chi_{2},23) = \textstyle \frac{1}{35}
(-6 \nu_{5a} + 20) \geq 0,
$$
that has no integral solutions such that all
$\mu_{i}(u,\chi_{j},23)$ are non-negative integers.

\noindent $\bullet$ Let $u$ be a unit of order $46$. We have
$\nu_{2a}+\nu_{23a}+\nu_{23b}=1$ by \eqref{E:1} and Proposition~\ref{P:4}. 
We have two cases: $\chi(u^2) \in \{ \chi(23a),
\chi(23b) \}$, but in both of them
$$
\mu_{0}(u,\chi_{9},*) = \textstyle \frac{1}{46} ( 286 \nu_{2a} +
266) \geq 0; \quad \mu_{23}(u,\chi_{9},*) = \textstyle
\frac{1}{46} (-286 \nu_{2a} + 240) \geq 0,
$$
that has no integral solutions such that all
$\mu_{i}(u,\chi_{j},*)$ are non-negative integers.

\noindent $\bullet$ Let $u$ be a unit of order $55$. We have
$\nu_{5a}+\nu_{11a}+\nu_{11b}=1$ by \eqref{E:1} and Proposition~\ref{P:4}. 
As before, we need to consider 20 cases given by part
(vii) of the Theorem, and all of them give us the same system of
inequalities
$$
\mu_{0}(u,\chi_{2},2) = \textstyle \frac{1}{55} (40 \nu_{5a} + 15)
\geq 0; \quad \mu_{5}(u,\chi_{2},2) = \textstyle \frac{1}{55} (-4
\nu_{5a} + 15) \geq 0,
$$
that has no integral solutions such that all
$\mu_{i}(u,\chi_{j},2)$ are non-negative integers.

\noindent $\bullet$ Let $u$ be a unit of order $69$. We have
$\nu_{3a}+\nu_{23a}+\nu_{23b}=1$ by \eqref{E:1} and Proposition~\ref{P:4}. 
We have two cases: $\chi(u^3) \in \{ \chi(23a),
\chi(23b) \}$, but both of them gives
\[
\begin{split}
\mu_{0}(u,\chi_{3},*) & = \textstyle \frac{1}{69} ( -44 (\nu_{23a}+\nu_{23b}) + 23) \geq 0;\\
\mu_{23}(u,\chi_{3},*) & = \textstyle \frac{1}{69} (22
(\nu_{23a}+\nu_{23b}) + 23) \geq 0,
\end{split}
\]
that has no integral solutions such that all
$\mu_{i}(u,\chi_{j},*)$ are non-negative integers.

\noindent $\bullet$ Let $u$ be a unit of order $77$. Clearly,
$\nu_{7a}+\nu_{7b}+\nu_{11a}+\nu_{11b}=1$ by \eqref{E:1} and
Proposition~\ref{P:4}. We need to consider 80 cases, determined by
parts (v) and (vii) of the Theorem. All of them give us the same
system of inequalities
\[
\begin{split}
\mu_{0}(u,\chi_{2},2) & = \textstyle \frac{1}{77} (30 (\nu_{7a}+\nu_{7b}) + 14) \geq 0; \\
\mu_{0}(u,\chi_{4},2) & = \textstyle \frac{1}{77} (-90
(\nu_{7a}+\nu_{7b}) +35) \geq 0,
\end{split}
\]
that has no integral solutions such that all
$\mu_{0}(u,\chi_{j},2)$ are non-negative integers.

\noindent $\bullet$ Let $u$ be a unit of order $115$. We have
$\nu_{5a}+\nu_{23a}+\nu_{23b}=1$ by \eqref{E:1} and Proposition~\ref{P:4}. 
Consider two cases: $\chi(u^5) \in \{ \chi(23a),
\chi(23b) \}$. Both of them give
\[
\begin{split}
\mu_{0}(u,\chi_{3},*) & = \textstyle \frac{1}{115} ( -88 (\nu_{23a}+\nu_{23b}) + 23) \geq 0;\\
\mu_{23}(u,\chi_{3},*) & = \textstyle \frac{1}{115} (22
(\nu_{23a}+\nu_{23b}) + 23) \geq 0,
\end{split}
\]
that has no integral solutions such that all
$\mu_{i}(u,\chi_{j},*)$ are non-negative integers.

\noindent $\bullet$ Let $u$ be a unit of order $161$. We have
$\nu_{7a}+\nu_{7b}+\nu_{23a}+\nu_{23b}=1$ by \eqref{E:1} and
Proposition~\ref{P:4}. Consider 8 cases, determined by possible
values of $\chi(u^{23})$ and $\chi(u^7)$. All of them give us the
same system of inequalities
\[
\begin{split}
\mu_{0}(u,\chi_{5},*) & = \textstyle \frac{1}{161} (-132 (\nu_{7a}+\nu_{7b}) + 224) \geq 0; \\
\mu_{0}(u,\chi_{9},*) & = \textstyle \frac{1}{161} (132
(\nu_{7a}+\nu_{7b}) +259) \geq 0,
\end{split}
\]
that has no integral solutions such that all
$\mu_{0}(u,\chi_{j},*)$ are non-negative integers.

\noindent $\bullet$ Let $u$ be a unit of order $253$. Again,
$\nu_{11a}+\nu_{11b}+\nu_{23a}+\nu_{23b}=1$ by \eqref{E:1} and
Proposition~\ref{P:4}. We need to consider 40 cases, determined by
possible values of $\chi(u^{23})$ and $\chi(u^{11})$. All of them
give us the same system of inequalities
\[
\begin{split}
\mu_{0}(u,\chi_{2},*)  & = \textstyle \frac{1}{253} (-220 (\nu_{23a}+\nu_{23b}) ) \geq 0; \\
\mu_{23}(u,\chi_{2},*) & = \textstyle \frac{1}{253} (  22 (\nu_{23a}+\nu_{23b}) ) \geq 0; \\
\mu_{11}(u,\chi_{6},*) & = \textstyle \frac{1}{253} ( -10 (\nu_{23a}+\nu_{23b}) + 230 ) \geq 0, \\
\end{split}
\]
that has no integral solutions such that all
$\mu_{0}(u,\chi_{j},*)$ are non-negative integers.

\subsection*{Acknowledgment}
The authors are grateful to Dr. Steve Linton and Dr. Tom Kelsey
from the University of St Andrews for their advice with
computational issues.

\bibliographystyle{alpha}
\bibliography{M23_final_Comm_Algebra}

\end{document}